\documentclass[12pt]{amsart}
\usepackage{geometry} 
\usepackage{graphicx}
\usepackage{amsmath}
\usepackage{amsfonts}
\usepackage{mathrsfs}
\usepackage{amssymb}
\usepackage{wrapfig}
\usepackage{bbold}
\usepackage{xcolor, soul}
\usepackage{mdframed}
\usepackage{enumitem}
\usepackage{hyperref}

\newtheorem{thm}{Theorem}
\newtheorem{lem}{Lemma}
\newtheorem{cor}{Corollary}
\newtheorem{prop}[cor]{Proposition}

\newtheorem{defn}{Definition}

\newenvironment{pf}{\noindent {\em Proof}.\ \ }{\hspace*{\fill}\rule{.5ex}{1.4ex}\,}

\numberwithin{equation}{section}
\numberwithin{defn}{section}

\newcommand{\real}{\mathbb{R}}
\newcommand{\complex}{\mathbb{C}}

\newcommand{\lloc}{L^1_{\mbox{\tiny loc}}(\mathbb{R})}
\newcommand{\llocgen}{L^1_{\mbox{\tiny loc}}}
\newcommand{\acloc}{AC_{\mbox{\tiny loc}}(\mathbb{R})}
\newcommand{\aclocgen}{AC_{\mbox{\tiny loc}}}
\newcommand{\atensorop}[1]{T^{(#1)}}
\newcommand{\atensor}[3]{T^{(#1)}_{#2,#3}}
\newcommand{\atensorjfg}{\atensor{j}{f}{g}}
\newcommand{\rsphere}{\widehat{\mathbb{C}}}
\DeclareMathOperator{\aut}{\rm Aut}
\newcommand{\autriem}{\aut\rsphere}
\DeclareMathOperator{\SL}{\rm SL}

\DeclareMathOperator{\GL}{\rm GL}
\DeclareMathOperator{\PSL}{\rm PSL}
\newcommand{\ster}{\mathcal{S}}
\newcommand{\X}{\mathbf{X}}
\newcommand{\C}{\mathbf{C}}
\newcommand{\M}{\mathbf{M}}
\newcommand{\Y}{\mathbf{Y}}
\DeclareMathOperator{\tr}{{\rm tr}}

\title{Solution of the scalar Riccati equation}

\author{Peter C.~Gibson}
\address{Dept.~of Mathematics \& Statistics\\
York University\\
Toronto, {ON}, M3J1P3, Canada}
\email{pcgibson@yorku.ca}

\date{August 4, 2025} 

\begin{document}

\maketitle

\begin{abstract}
The scalar Riccati equation is a prototypical nonlinear ODE having diverse mathematical connections. In the centuries since its initial formulation, a standard textbook theory has emerged according to which the general solution may be determined if a particular solution is known; but no general method exists to determine a particular solution explicitly, except in sporadic special cases. 
The purpose of the present article is to solve the scalar Riccati equation in general form, as well as the general linear ODE of second order, directly by explicit construction.  
In the case of the Riccati equation, the solution 
sets up a bijective correspondence between triples of locally integrable functions on the real line,
and locally absolutely continuous paths through the identity in the automorphism group of the Riemann sphere.  
As applications of the results, we obtain an explicit solution to the one-dimensional Schr\"odinger equation, an inversion formula for the Miura transform, and a new formula for Airy functions. 
\end{abstract}

\vspace*{12pt}


%
\tableofcontents


\section{Introduction\label{sec-introduction}}

We solve explicitly the scalar Riccati equation
\begin{equation}\label{riccati}
y^\prime=fy^2+gy+h
\end{equation}
where $f,g,h$ are arbitrary locally integrable complex valued functions on $\real$.  Methods to integrate (\ref{riccati}) have heretofore been established only in sporadic special cases. 
The technical key to our solution is a newly defined special operator, called the bivariate exponential, that generalizes the (univariate) exponential primitive operator 
\begin{equation}\nonumber
\psi\mapsto e^{\int_0^x\psi}=1+\sum_{j=1}^\infty\thickspace\int\limits_{0<s_1<\cdots<s_j<x}\psi(s_1)\cdots\psi(s_j)\,ds_1\cdots ds_j
\end{equation}
 in a natural way, interleaving two functions in place of one:
 \begin{equation}\nonumber
 \begin{split}
 &\mbox{\tiny $j$ terms, alternating $\psi$ with $\eta$ }\\
 (\psi,\eta)\mapsto E_{\psi,\eta}(x)=1+\sum_{j=1}^\infty\thickspace \int\limits_{0<s_1<\cdots<s_j<x}&\overbrace{\psi(s_1)\eta(s_1)\psi(s_3)\cdots}\,ds_1\cdots ds_j.
 \end{split}
 \end{equation}
The solution to (\ref{riccati}) with arbitrary initial value $y(0)=y_0\in\complex$ has a surprisingly
simple formulation in terms of the above operator, as follows.
Setting
\begin{equation}\nonumber
C_{\psi,\eta}=\frac{1}{2}\bigl(E_{\psi,\eta}+E_{-\psi,-\eta}\bigr)\quad\mbox{ and }\quad S_{\psi,\eta}=\frac{1}{2}\bigl(E_{\psi,\eta}-E_{-\psi,-\eta}\bigr),
\end{equation}
\begin{equation}\nonumber
y(x)=e^{\int_0^x\negthickspace g}\thinspace\frac{y_0C_{F,H}(x)+S_{H,F}(x)}{y_0S_{F,H}(x)+C_{H,F}(x)}\quad\mbox{ where }\quad F=-fe^{\int_0^xg}\mbox{ and }H=he^{-\int_0^xg}.
\end{equation}
Minor modifications extend the above formulas to all $x\in\real$, and the corresponding solution is global for almost every $y_0\in\complex$. 

By the same token, for any locally integrable $\alpha,\beta$, the general solution to
\begin{equation}\nonumber
y^{\prime\prime}+\alpha y^\prime+\beta y=0
\end{equation}
is expressed  very simply in terms of the operators $C$ and $S$ by the formula
\begin{equation}\nonumber
y=y(0)C_{f,g}+y^\prime(0)S_{g,f}\quad\mbox{ where }\quad f=-\beta e^{\int_0^x\alpha}\mbox{ and }g=e^{-\int_0^x\alpha}. 
\end{equation}
Particular applications of these results yield a new formula for Airy functions, an explicit solution to the one-dimensional Schr\"odinger equation, and an inversion formula for the Miura transform.  

\subsection{Background, motivation and goals\label{sec-background}}
The Riccati equation (\ref{riccati}) has a long history dating at least to 1720  \cite{Bi:1996}, with early contributions to its analysis being made by Euler, Abel and Liouville \cite{Eu:1785,Eu:1818,Kh:2006,Ab:1881,Li:1841}.  
Its recent importance stems from diverse connections, including to calculus of variations, control theory, algebraic geometry and mathematical physics \cite{Re:1972,Ze:2000,LlVa:2014,GaArTo:2016}.  The references in the latter works, as well as those in \cite{Bi:1996}, give an indication of the considerable scope of literature relating to (\ref{riccati}). 
The present article was motivated by a desire to invert the Miura transform $\alpha\mapsto\alpha^2+\alpha^\prime$, which arises in connection to nonlinear waves \cite{Mi:1968,AbKaNeSe:1974,KaPeShTo:2005}. The crucial technical insight underpinning our solution stems from recent analysis of the Schr\"odinger equation \cite{Gi:SIMA2024}.  

Irrespective of context, be it control theory or mathematical physics, the most basic problem inherent in the Riccati equation is to solve it, or to prove no solution exists.  As evidenced for example by \cite{ZwDo:2022,PoZa:2003,En:2025}, established theory concerning the basic problem has not substantially changed since the original publication of Kamke's compendium \cite{Ka:1977}.  To summarize briefly, one part of established theory consists of the general facts that: one can determine the general solution to (\ref{riccati}) given a particular solution; the cross ratio of any four distinct particular solutions is constant; equation (\ref{riccati}) may be transformed into a linear second order equation, or into a canonical form in which $g=0$.  A second part of established theory consists of a list of sporadic special cases in which a particular solution to (\ref{riccati}) may be determined explicitly by elementary combinations of quadratures.  Some of these sporadic special cases are quite complicated, and every so often a new one not listed by Kamke is discovered, such as \cite{St:1991,MaHa:2013}.  Liouville determined precisely when solution by quadratures is possible for a special subclass of Riccati equations \cite{Li:1841}.
 But it turns out that from the point of view of the basic problem these sporadic special cases are something of a red herring. 

A main goal of the present paper is to show that one can write down the general solution to (\ref{riccati}) directly, without recourse to complicated transformations, or to a taxonomy of special subclasses.  Indeed, as described above, solution of the Riccati equation is surprisingly straightforward given the right tool, namely the bivariate exponential operator.  A second goal is to place solution of the Riccati equation in its natural mathematical setting, wherein continuous paths through the identity in the automorphism group of the Riemann sphere generate solutions to a corresponding Riccati equation.  More precisely, there is a bijective ``Riccati correspondence" between triples $(f,g,h)$ of locally integrable coefficients, and continuous paths through the identity in $\aut\rsphere$ that generate the solutions to (\ref{riccati}) for arbitrary initial data in $\complex$.  We construct this correspondence explicitly.

\subsection{Organization of the paper\label{sec-organization}}

Section~\ref{sec-preliminaries}  lays out technical preliminaries, notation and definitions that are needed in the rest of the paper. In particular, the notions of standard solution, projective solution and singular curve are made precise.

The paper's main results are stated  in \S\ref{sec-main}. The bivariate exponential operator is defined in technically precise terms in \S\ref{sec-bivariate-exponential}, with Lemma~\ref{lem-differential-formula}  encapsulating its crucial properties. Theorem~\ref{thm-matrix-equation} in \S\ref{sec-companion-solution} gives the solution to the $2\times2$ matrix equation associated with (\ref{riccati}).  Theorem~\ref{thm-linear} in \S\ref{sec-linear} solves the general homogeneous linear equation of second order, and Theorem~\ref{thm-riccati-solution} in \S\ref{sec-riccati-solution} solves (\ref{riccati}) itself.  The Riccati correspondence is spelled out in \S\ref{sec-general-properties}.  

Detailed proofs of the main results are given in \S\ref{sec-proofs}.  

Section~\ref{sec-examples} applies the main theorems to several specific cases,
both to illustrate their utility, and to obtain results of independent interest. A solution originally discovered by Abel, valid for a special subclass of Riccati equations, is recovered (and slightly generalized) in \S\ref{sec-abel}. 
In \S\ref{sec-Airy} we use Theorem~\ref{thm-matrix-equation} to obtain a new formula for Airy functions.  
In \S\ref{sec-schroedinger} we use Theorem~\ref{thm-matrix-equation} to solve the one-dimensional Schr\"odinger equation explicitly, thereby replicating an earlier result which was originally obtained by a much more complicated and technically difficult chain of reasoning.  In \S\ref{sec-miura}, we use Theorem~\ref{thm-riccati-solution} to derive an explicit inverse formula for the Miura transform, which is a new result. 

Section~\ref{sec-remarks} rounds out the paper with some technical remarks and general discussion.

\section{Technical preliminaries\label{sec-preliminaries}}

\subsection{Automorphisms of the Riemann sphere\label{sec-automorphisms}}

A useful elementary fact from classical analysis is that any invertible $2\times 2$ matrix 
\begin{equation}\label{gl2c}
M=\begin{pmatrix}a&b\\ c&d\end{pmatrix}\in\GL(2,\complex)
\end{equation}
determines an automorphism of the Riemann sphere $\rsphere=\complex\cup\{\infty\}$ via the formula
\begin{equation}\label{lft}
\varphi_M(z)=\frac{az+b}{cz+d}\qquad(z\in\rsphere).
\end{equation}
And every member of $\autriem$ arises this way.
Two matrices $M_1,M_2\in\GL(2,\complex)$ determine the same automorphism $\varphi_{M_1}=\varphi_{M_2}$ if and only if they are nonzero scalar multiples of one another. 
Thus $\autriem$ is represented by $\PSL(2,\complex)$, elements of which are pairs $\{\pm M\}$, where $M\in\SL(2,\complex)$, the group of $2\times 2$ matrices of determinant 1. 

A continuous matrix valued function $M:\real\rightarrow\SL(2,\complex)$ thus determines a path 
\begin{equation}\label{path}
\varphi_M:\real\rightarrow\autriem.
\end{equation}
For each fixed $z\in\rsphere$, the path (\ref{path}) in turn determines a function
\begin{equation}\label{solution-function}
\varphi_M(z):\real\rightarrow\rsphere,\quad x\mapsto \varphi_{M(x)}(z).  
\end{equation} 
Solutions to the Riccati equation are known to be functions of the form (\ref{solution-function}).  But there is no known method to construct $M$ explicitly.  That is a principal objective of the present paper. 

Unless stated otherwise, it is to be understood in the present article that the Riemann sphere carries the metric from $\real^3$ induced by stereographic projection
\begin{equation}\label{stereographic}
\ster:\rsphere\rightarrow S^2\subset\real^3,\qquad \ster(z)=
\left\{
\begin{array}{cc}\frac{1}{|z|^2+1}\bigl(2\Re z,2\Im z,|z|^2-1\bigr)&\mbox{ if }z\in\complex\\
(0,0,1)&\mbox{ if }z=\infty
\end{array}
\right..
\end{equation}

\subsection{Notation and terminology\label{sec-notation}}
Denote by $\lloc$ the space of locally integrable complex valued functions on the real line, and by $\acloc$ the space of locally absolutely continuous complex valued functions. More generally, given a subset $Y$ either of $\complex^{2\times 2}$ or of $\real^3$, the set of locally integrable functions from $\real$ to $Y$ is denoted $\llocgen(\real,Y)$, and the set of locally absolutely continuous functions from $\real$ to $Y$ is denoted $\aclocgen(\real,Y)$.  Let 
\begin{equation}\nonumber
I=\begin{pmatrix}1&0\\ 0&1\end{pmatrix}
\end{equation}
denote the $2\times2$ identity matrix. 
Denote
\begin{equation}\label{spaces}
\begin{split}
\X&=\lloc\times\lloc\times\lloc\\
\C&=\llocgen\bigl(\real,\mathfrak{sl}(2,\complex)\bigr)\\
\M&=\left\{M\in\aclocgen(\real,\SL(2,\complex)\,\left|\,M(0)=I\right.\right\}\\
\Y&=\left\{\left.\varphi_M:\real\rightarrow\autriem\,\right|\,M\in\M\right\}
\end{split}
\end{equation}
Note that since $M(0)=I$ for any $M\in\M$, it follows that for any $\varphi_M\in\Y$, 
\begin{equation}\label{identity}
\varphi_{M(0)}(z)=z\qquad(z\in\rsphere). 
\end{equation}


By a \emph{standard solution} to equation (\ref{riccati}), we mean a function $y\in\acloc$ for which (\ref{riccati}) holds almost everywhere.  It follows from the classical theory due to Carath\'eodory that there is at most one standard solution to (\ref{riccati}) consistent with any given initial condition of the form $y(0)=y_0$, where $y_0\in\complex$. See \cite[Thm.~I.5.3]{Ha:1980}. A standard solution $y$ cannot take the value $y(x)=\infty$, since $y$ is necessarily discontinuous at such a point $x$.  To accommodate the latter possibility, we define a second, slightly weaker notion of solution based on stereographic projection (\ref{stereographic}), as follows. 

A function $y:\real\rightarrow\rsphere$ will be referred to as a \emph{projective solution} to (\ref{riccati}) if 
\begin{equation}\label{condition}
\ster y\in\aclocgen(\real,S^2)
\end{equation}
and $(\ref{riccati})$ holds almost everywhere. Thus every standard solution is a projective solution, but a projective solution may take the value $\infty$, whereas a standard solution cannot.  

Referring to (\ref{spaces}), an element $\varphi_M\in\Y$ will be said to \emph{generate solutions} to (\ref{riccati}), or to be a \emph{solution generator} for (\ref{riccati}), if for almost every $y_0\in\complex$, the function 
\begin{equation}\label{solution-generator}
y(x)=\varphi_{M(x)}(y_0)\qquad(x\in\real)
\end{equation}
is a standard solution to (\ref{riccati}).  Observe that the function (\ref{solution-generator}) satisfies the initial condition $y(0)=y_0$ by construction.  It follows from the Carath\'eodory uniqueness theorem that a Riccati equation can have at most one solution generator, since an element of $\autriem$ is uniquely determined by its values on any three distinct points.  In detail, let $\varphi_{M_1},\varphi_{M_2}\in\Y$ generate solutions to (\ref{riccati}). Then there exists a set $S\subset\complex$ of measure zero such that for every $z\in\complex\setminus S$, 
\begin{equation}\nonumber
y_1=\varphi_{M_1}(z)\quad\mbox{ and }\quad y_2=\varphi_{M_2}(z)
\end{equation}
are standard solutions to (\ref{riccati}). Since $y_1(0)=y_2(0)=z$, it follows by uniqueness of standard solutions that $y_1(x)=y_2(x)$ for every $x\in\real$.  Thus 
\begin{equation}\label{generator-identity}
\varphi_{M_1(x)}(z)=\varphi_{M_2(x)}(z)
\end{equation}
for every $z\in\complex\setminus S$ and $x\in\real$.  Since for each fixed $x\in\real$, (\ref{generator-identity}) holds for more than three distinct $z\in\complex$, it follows that 
\begin{equation}\label{automorphism-identity}
\varphi_{M_1(x)}=\varphi_{M_2(x)}.
\end{equation}
Since (\ref{automorphism-identity}) holds for every $x\in\real$, it follows in turn that $\varphi_{M_1}=\varphi_{M_2}$, as claimed. 

\subsection{The singular curve\label{sec-singular-curve}} 

Define the \emph{singular curve} $\Sigma_M\subset\complex$ associated to a matrix function
\begin{equation}\label{matrix-function-2}
M=\begin{pmatrix}a&b\\ c&d\end{pmatrix}\in\M
\end{equation}
to be the set
\begin{equation}\label{singular-curve-0}
\Sigma_M=\left\{\left.-d(x)/c(x)\,\right|\,x\in\real\mbox{ and }c(x)\neq0\right\}.
\end{equation}

\begin{prop}\label{prop-measure-0}
For every $M\in\M$, the singular curve $\Sigma_M$ has measure zero in $\complex$. 
\end{prop}
\begin{pf}
Observe that the set of all $x\in\real$ such that $c(x)\neq 0$ is a countable union of open intervals, each of which is a countable union of (connected) compact intervals.  By absolute continuity, the set of points $-d(x)/c(x)$ over such a compact interval is a rectifiable curve, and hence has measure zero. Since $\Sigma_M$ is covered by a countable collection of such sets, $\Sigma_M$ itself has measure zero. 
\end{pf}

\section{Main Results\label{sec-main}}

\subsection{Formal definition of the bivariate exponential operator\label{sec-bivariate-exponential}}

 Let $P$ denote
\begin{equation}\label{primitive}
P:\lloc\rightarrow\acloc,\qquad Pf(x)=\int_0^xf.
\end{equation}
Thus $Pf$ is the primitive of $f$, normalized such that $Pf(0)=0$.  Solution of Riccati equation involves the following 
generalization of the exponential primitive operator $f\mapsto e^{Pf}$. 

For each integer $j\geq 1$, define
\begin{equation}\label{Tj}
\begin{split}
\atensorop{j}:\lloc\times\lloc\rightarrow L^1_{\mbox{\tiny loc}}(\real^j),\qquad\atensorjfg=&\underbrace{f\otimes g\otimes f\otimes\cdots\thickspace}\\
&\mbox{\tiny $j$ terms, alternating $f$ with $g$ }
\end{split}
\end{equation}
Thus, for example, $\atensor{1}{f}{g}(s_1)=f(s_1)$, $\atensor{2}{f}{g}(s_1,s_2)=f(s_1)g(s_2)$, 
\begin{equation}\nonumber
\atensor{3}{f}{g}(s_1,s_2,s_3)=f(s_1)g(s_2)f(s_3),\quad\atensor{4}{f}{g}(s_1,s_2,s_3,s_4)=f(s_1)g(s_2)f(s_3)g(s_4),\mbox{ etc. }
\end{equation}
For $x\in\real$ and $f,g\in\lloc$, set 
\begin{equation}\label{Aj}
A_j(f,g,x)=\left\{
\begin{array}{cc}
\displaystyle\int\limits_{0<s_1<\cdots<s_j<x}\atensorjfg(s_1,\ldots,s_j)\,ds_1\cdots ds_j&\mbox{ if }x\geq 0\\
\rule{0pt}{21pt}(-1)^j\displaystyle\int\limits_{x<s_j<\cdots<s_1<0}\atensorjfg(s_1,\ldots,s_j)\,ds_1\cdots ds_j&\mbox{ if }x\leq 0
\end{array}
\right.
\end{equation}
and set $A_0(f,g,x)=1$. 

\begin{defn}\label{defn-bivariate-exponential}
Define the \emph{bivariate exponential} operator
\begin{equation}\label{epo1}
E:\lloc\times\lloc\rightarrow\acloc,\qquad (f,g)\mapsto E_{f,g}
\end{equation}
by the formula
\begin{equation}\label{epo}
E_{f,g}(x)=1+\sum_{j=1}^\infty A_j(f,g,x)\qquad(x\in\real)
\end{equation}
and set
\begin{equation}\label{cs}
C_{f,g}=\frac{1}{2}\left(E_{f,g}+E_{-f,-g}\right),\qquad S_{f,g}=\frac{1}{2}\left(E_{f,g}-E_{-f,-g}\right).
\end{equation} 
\end{defn}
Observe that $E_{f,g}(0)=C_{f,g}(0)=1$ and $S_{f,g}(0)=0$, irrespective of $f,g\in\lloc$. 
That the bivariate exponential has codomain $\acloc$ requires proof.  Indeed,
\begin{lem}
\label{lem-differential-formula}
For every $f,g\in\lloc$, the functions $E_{f,g}, C_{f,g}$ and $S_{f,g}$ are locally absolutely continuous, with
\begin{equation}\label{differential-formula-cs}
\left(C_{f,g}\right)^\prime=gS_{f,g}\quad\mbox{ and }\quad \left(S_{f,g}\right)^\prime=fC_{f,g}.
\end{equation}
\end{lem}
The differential formulas (\ref{differential-formula-cs}) are key to solving the general Riccati equation in explicit form. Lemma~\ref{lem-differential-formula} is proved in \S\ref{sec-lemma-2-proof}. 

It follows from Lemma~\ref{lem-differential-formula} that for any $f\in\lloc$,
\begin{equation}\label{diagonal-case}
E_{f,f}^\prime=C_{f,f}^\prime+S_{f,f}^\prime=fS_{f,f}+fC_{f,f}=fE_{f,f}. 
\end{equation}
And $E_{f,f}(0)=1$. Therefore, 
\begin{equation}\label{diagonal-values}
E_{f,f}=e^{Pf},\qquad C_{f,f}=\cosh Pf,\quad\mbox{ and }\quad S_{f,f}=\sinh Pf. 
\end{equation}
Thus the bivariate exponential generalizes the exponential primitive operator
\begin{equation}\label{exponential-primitive}
\lloc\ni f\mapsto e^{Pf}\in\acloc,
\end{equation}
as claimed.
In relation to the bivariate exponential $E$, the operators $C$ and $S$ are analogous to hyperbolic cosine and hyperbolic sine, respectively.

\subsection{Solution of the associated matrix equation\label{sec-companion-solution}}

Define the \emph{solution matrix} operator $\psi$ on $\X$
by the formula
\begin{equation}\label{psi}
\psi(f,g,h)=\begin{pmatrix}e^{\frac{1}{2}Pg}C_{F,H}&e^{\frac{1}{2}Pg}S_{H,F}\\
e^{-\frac{1}{2}Pg}S_{F,H}&e^{-\frac{1}{2}Pg}C_{H,F}\end{pmatrix}\quad\mbox{ where }\quad F=-e^{Pg}f \mbox{ and } H=e^{-Pg}h.
\end{equation}

\begin{thm}\label{thm-matrix-equation}
Let $(f,g,h)\in\X$ be arbitrary, and set $M=\psi(f,g,h)$. 
Then $M\in\M$; i.e., each entry of $M$ is locally absolutely continuous, and $\det M\equiv 1$.  
Moreover,
\begin{equation}\label{matrix-equation-2}
M^\prime=\mathfrak{m}M\quad\mbox{ and }\quad M(0)=I,\quad\mbox{ where }\quad\mathfrak{m}=\begin{pmatrix}\frac{1}{2}g&h\\ -f&-\frac{1}{2}g\end{pmatrix}.
\end{equation}
\end{thm}

\subsection{Solution of homogenous linear second order ODEs\label{sec-linear}}
Given an arbitrary complex $2\times 2$ matrix $\mathfrak{n}$, set $\mathfrak{m}=\mathfrak{n}-\frac{1}{2}\tr\mathfrak{n}I$.  If $M$ satisfies the matrix equation $M^\prime=\mathfrak{m}M$, then $N=e^{\frac{1}{2}P\tr\mathfrak{n}}M$ satisfies $N^\prime=\mathfrak{n}N$.  In this way, Theorem~\ref{thm-matrix-equation} extends to arbitrary coefficient matrices $\mathfrak{n}$.  In particular, any homogeneous linear second order ODE can be represented as a $2\times 2$ system whose general solution can be obtained using Theorem~\ref{thm-matrix-equation}, as follows.  
\begin{thm}\label{thm-linear}
For any $\alpha,\beta\in\lloc$, the differential equation 
\begin{equation}\label{linear}
y^{\prime\prime}+\alpha y^\prime+\beta y=0
\end{equation}
has general solution
\begin{equation}\label{linear-solution}
y=y(0)C_{f,g}+y^\prime(0)S_{g,f}\quad\mbox{ where }\quad f=-\beta e^{P\alpha}\mbox{ and }g=e^{-P\alpha}.
\end{equation}
\end{thm}
Note that formula (\ref{linear-solution}) can be verified directly using Lemma~\ref{lem-differential-formula}. 

\subsection{Solution of the Riccati equation\label{sec-riccati-solution}}

As a consequence of Theorem~\ref{thm-matrix-equation}, every Riccati equation (\ref{riccati}) has a solution generator that can be formulated explicitly in terms of the operators $C$ and $S$ of Definition~\ref{defn-bivariate-exponential}.  Under very mild conditions on the coefficient $f$ of $y^2$, the solution generator produces a global projective solution for every initial condition, without exception. 
\begin{thm}\label{thm-riccati-solution}
Let $(f,g,h)\in\X$ be arbitrary, and set $M=\psi(f,g,h)$ as per (\ref{psi}). Then
for every $y_0\in\complex\setminus\Sigma_M$, the function 
\begin{equation}\label{global-solution}
y(x)=\varphi_{M(x)}(y_0)\qquad(x\in\real)
\end{equation}
is the unique solution in the standard sense to the Riccati equation
\begin{equation}\label{riccati-2}
y^\prime=fy^2+gy+h,\qquad y(0)=y_0.
\end{equation}
For every  $y_0\in\Sigma_M$, the function (\ref{global-solution}) satisfies $\ster y\in\aclocgen(\real,S^2)$. 
If $f^{-1}(0)$ has measure zero, then for every $y_0\in\Sigma_M$, equation (\ref{riccati-2}) holds almost everywhere, and the function $y$ defined by (\ref{global-solution}) is the unique projective solution to (\ref{riccati-2}). 
\end{thm}

Written out more fully, the solution (\ref{global-solution}) to the Riccati equation (\ref{riccati-2}) is 
\begin{equation}\label{y-full}
y=e^{Pg}\frac{y_0C_{F,H}+S_{H,F}}{y_0S_{F,H}+C_{H,F}}\quad\mbox{ where }\quad F=-e^{Pg}f\mbox{ and }H=e^{-Pg}h. 
\end{equation}

\subsection{The Riccati correspondence\label{sec-general-properties}}

Referring to (\ref{spaces}), there is a bijective correspondence $\Phi:\X\rightarrow\Y$ mapping the coefficients of a Riccati equation to its unique solution generator.  
With Theorems~\ref{thm-matrix-equation} and \ref{thm-riccati-solution} in hand, this \emph{Riccati correspondence} can be built up as a composition
\begin{equation}\label{composition}
\X\rightarrow\C\rightarrow\M\rightarrow\Y
\end{equation}
as follows. 

Given $(f,g,h)\in\X$, set 
\begin{equation}\label{m}
\mathfrak{m}=\begin{pmatrix}\frac{1}{2}g&h\\ -f&-\frac{1}{2}g\end{pmatrix}\in\C. 
\end{equation}
Consider the $2\times2$ matrix equation
\begin{equation}\label{matrix-equation-1}
M^\prime=\mathfrak{m}M,\qquad M(0)=I.
\end{equation}
Theorem~\ref{thm-matrix-equation} provides a continuous global solution to (\ref{matrix-equation-1})
 for arbitrary $\mathfrak{m}\in\C$. Standard theory implies this global solution is unique.
Again by Theorem~\ref{thm-matrix-equation}, the unique solution $M$ to (\ref{matrix-equation-1}) belongs to $\M$.  Hence $\varphi_M\in\Y$, since $\mathfrak{m}\in\C$ has trace zero.
Theorem~\ref{thm-riccati-solution} implies $\varphi_M$ is the solution generator for (\ref{riccati}).  This can be seen directly, as follows. To streamline notation, denote 
\begin{equation}\label{matrix-function-3}
M=\begin{pmatrix}a&b\\ c&d\end{pmatrix}
\end{equation}
so that, by (\ref{matrix-equation-1}),
\begin{equation}\label{solution-formulas}
\begin{pmatrix}a^\prime&b^\prime\\ c^\prime&d^\prime\end{pmatrix}=
\begin{pmatrix}\frac{1}{2}ga+hc&\frac{1}{2}gb+hd\\ -fa-\frac{1}{2}gc&-fb-\frac{1}{2}gd\end{pmatrix}.
\end{equation}
Fix $y_0\in\complex$.  Differentiation of 
\begin{equation}\label{lft-2}
y=\varphi_M(y_0)=\frac{y_0a+b}{y_0c+d}
\end{equation}
with respect to the independent variable $x\in\real$, and application of the identities (\ref{solution-formulas}), easily yields (\ref{riccati}) at points $x$ where 
\begin{equation}\label{condition-2}
y_0c(x)+d(x)\neq0. 
\end{equation}
Observe that (\ref{condition-2}) holds for all $x\in\real$, unless $y_0$ lies on the singular curve $\Sigma_M$.  And by Proposition~\ref{prop-measure-0}, the singular curve has measure zero in $\complex$. 
Thus $\varphi_M$ is a solution generator for (\ref{riccati}), as claimed.

Inversely, let $\varphi_M\in\Y$ be arbitrary. Then $\varphi_M$ determines a unique matrix function
\begin{equation}\label{matrix-function-4}
M=\begin{pmatrix}a&b\\ c&d\end{pmatrix}\in\M.
\end{equation}
Set 
\begin{equation}\label{m-formula}
\mathfrak{m}=M^\prime M^{-1}=\begin{pmatrix}a^\prime d-b^\prime c&ab^\prime-a^\prime b\\ c^\prime d-cd^\prime&ad^\prime-bc^\prime\end{pmatrix}
\end{equation}
so that (\ref{matrix-equation-1}) holds.  Note that $\tr\mathfrak{m}(x)=0$ a.e., since $ad-bc=\det M\equiv 1$.
Setting
\begin{equation}\label{fgh-formulas}
f=cd^\prime-c^\prime d,\quad g=2(a^\prime d-b^\prime c)\quad\mbox{ and }\quad h=ab^\prime-a^\prime b,
\end{equation}
it follows that $(f,g,h)\in\X$ and $\Phi(f,g,h)=\varphi_M$. 
Thus the Riccati correspondence $\Phi:\X\rightarrow\Y$ described above is bijective.

\section{Analysis and proofs\label{sec-proofs}}

\subsection{Functions invariant under permutation\label{sec-invariant}}

The following technical fact is useful in analyzing integrals over simplices. 
\begin{prop}\label{prop-permutations}
Let $\sigma_j$ denote the symmetric group of permutations on $\{1,2,\ldots,j\}$, and let $F\in L^1_{\mbox{\tiny loc}}(\mathbb{R}^j)$.
If $F$ is invariant under permutations in the sense that
\begin{equation}\label{invariance}
F(s_1,\ldots,s_j)=F(s_{\pi(1)},\ldots,s_{\pi(j)})\quad\mbox{ for every }\quad (s_1,\ldots,s_j)\in\real^j\mbox{ and }\pi\in\sigma_j
\end{equation}
then for every $x_1<x_2$,
\begin{equation}\label{two-integrals}
\int_{x_1<s_1<\cdots<s_j<x_2}F(s_1,\ldots,s_j)\,ds_1\cdots ds_j=\frac{1}{j!}\int_{(x_1,x_2)^j}F(s_1,\ldots,s_j)ds_1\cdots ds_j. 
\end{equation}
\end{prop}
\begin{pf}
The $j!$ simplices consisting of all $(s_1,\ldots,s_j)\in\real^j$ such that 
\[
x_1<s_{\pi(1)}<\cdots<s_{\pi(j)}<x_2\quad\mbox{ where }\quad \pi\in\sigma_j
\]
are congruent to one another and cover the $j$-cube $(x_1,x_2)^j$ up to a set of measure zero.  Therefore
\begin{equation}\label{partition}
\begin{split}
\int_{(x_1,x_2)^j}F(s_1,\ldots,s_j)&ds_1\cdots ds_j\\
&=\sum_{\pi\in\sigma_j}\int_{x_1<s_{\pi(1)}<\cdots<s_{\pi(j)}<x_2}F(s_1,\ldots,s_j)\,ds_1\cdots ds_j\\
&=j!\int_{x_1<s_1<\cdots<s_j<x_2}F(s_1,\ldots,s_j)\,ds_1\cdots ds_j
\end{split}
\end{equation}
the latter equality by invariance of $F$.
\end{pf}

\subsection{Proof of Lemma~\ref{lem-differential-formula}\label{sec-lemma-2-proof}}

Referring to (\ref{Aj}), observe
\begin{equation}\label{Tj-bound}
\left|T_{f,g}(s_1,\ldots,s_j)\right|\leq T_{|f|+|g|,|f|+|g|}(s_1,\ldots,s_j). 
\end{equation}
Therefore Proposition~\ref{prop-permutations} applied to integration of the right-hand side of (\ref{Tj-bound}) yields
\begin{equation}\label{Aj-bound}
\left|A_j(f,g,x)\right|\leq \frac{1}{j!}\left(P(|f|+|g|)(x)\right)^j. 
\end{equation}
For $j\geq 1$, set 
\begin{equation}\label{Bj}
B_j(f,g,x)=\left\{
\begin{array}{cc}
f(x)A_{j-1}(f,g,x)&\mbox{ if $j$ is odd }\\
g(x)A_{j-1}(f,g,x)&\mbox{ if $j$ is even }
\end{array}\right.
\end{equation}
so that by (\ref{Aj-bound}),
\begin{equation}\label{Bj-bound}
\left|B_j(f,g,x)\right|\leq\frac{|f|+|g|}{(j-1)!}\left(P(|f|+|g|)(x)\right)^{j-1}. 
\end{equation}
Note that formula (\ref{Aj}) implies
\begin{equation}\label{Aj-integral-form}
A_j(f,g,x)=\int_0^xB_j(f,g,y)\,dy
\end{equation}
for every $x\in\real$ and $j\geq1$. Now, it follows from (\ref{Bj-bound}) that the series $\sum B_j(f,g,x)$ converges absolutely, pointwise almost everywhere, and that partial sums are collectively dominated by a locally integrable function:
\begin{equation}\label{partial-sum-bound}
\left|\sum_{j=1}^nB_j(f,g,x)\right|\leq\bigl(|f(x)|+|g(x)|\bigr)e^{P(|f|+|g|)(x)}\qquad(n\geq1).
\end{equation}
Dominated convergence implies the series $\sum B_j(f,g,x)$ is locally integrable, and 
\begin{equation}\label{E-AC}
\int_{0}^x\sum_{j=1}^\infty B_j(f,g,y)\,dy=\sum_{j=1}^\infty\int_{0}^xB_j(f,g,y)\,dy=\sum_{j=1}^\infty A_j(f,g,x)=E_{f,g}(x)-1. 
\end{equation}
Thus $E_{f,g}-1$ is the primitive of a locally integrable function, and $E_{f,g}$ is therefore locally absolutely continuous.  It follows that both $C_{f,g}$ and $S_{f,g}$ are also locally absolutely continuous. 

The bound (\ref{Aj-bound}) implies that the series defining $E_{f,g}$ converges absolutely, uniformly on compact sets.  The order of terms in the respective series defining $C_{f,g}$ and $S_{f,g}$ may therefore be rearranged as desired.  The resulting cancellations, together with dominated convergence, lead easily to the formulas
\begin{equation}\label{C-integral-formulas}
\begin{split}
C_{f,g}(x)&=1+\sum_{j=1}^\infty A_{2j}(f,g,x)\\
&=1+\sum_{j=1}^\infty\int_{0}^xg(y)A_{2j-1}(f,g,y)\,dy\\
&=1+\int_0^x\sum_{j=1}^\infty g(y)A_{2j-1}(f,g,y)\,dy
\end{split}
\end{equation}
and
\begin{equation}\label{S-integral-formulas}
\begin{split}
S_{f,g}(x)&=\sum_{j=1}^\infty A_{2j-1}(f,g,x)\\
&=\sum_{j=1}^\infty\int_{0}^xf(y)A_{2j-2}(f,g,y)\,dy\\
&=\int_0^x\sum_{j=1}^\infty f(y)A_{2j-2}(f,g,y)\,dy.
\end{split}
\end{equation}
Differentiation of (\ref{C-integral-formulas}) and (\ref{S-integral-formulas}) yields 
\begin{equation}\label{desired-formulas}
\frac{d}{dx}C_{f,g}(x)=g(x)S_{f,g}(x)\quad\mbox{ and }\quad\frac{d}{dx}S_{f,g}(x)=f(x)C_{f,g}(x).
\end{equation}
This completes the proof of Lemma~\ref{lem-differential-formula}. 

\subsection{Proof of Theorm~\ref{thm-matrix-equation}\label{sec-thm-1-proof}}

It follows from Lemma~\ref{lem-differential-formula} that each of the entries of $M=\psi(f,g,h)$ is locally absolutely continuous, and that $M$ satisfies (\ref{matrix-equation-2}).  The condition $M(0)=I$ is also straightforwardly verified from the definition of operators $C$ and $S$. It follows in turn, since $\mathfrak{m}$ as defined in (\ref{matrix-equation-2}) has trace 0, that $\det M$ is constant, the constant value being $\det M(0)=1$.  To see in detail that $\det M$ is constant, let $a,b,c,d\in\acloc$ and $\alpha,\beta,\gamma\in\lloc$.  If matrices
\begin{equation}\label{M-with-entries}
M=\begin{pmatrix}a&b\\ c&d\end{pmatrix}\quad\mbox{ and }\quad\mathfrak{m}=\begin{pmatrix}\beta&\gamma\\ \alpha&-\beta\end{pmatrix}
\end{equation}
satisfy the identity $M^\prime=\mathfrak{m}M$ a.e., then application of this identity to 
\begin{equation}\label{det-expansion}
\bigl(\det M\bigr)^\prime=a^\prime d+ad^\prime-b^\prime c-bc^\prime
\end{equation}
yields that $\bigl(\det M\bigr)^\prime=0$ almost everywhere. 
Therefore $\det M$ is constant, as claimed. This completes the proof of Theorem~\ref{thm-matrix-equation}.

As mentioned already, Theorem~\ref{thm-linear} may be verified directly using Lemma~\ref{lem-differential-formula}.  Thus we proceed to Theorem~\ref{thm-riccati-solution}.

\subsection{Proof of Theorem~\ref{thm-riccati-solution}\label{sec-thm-2-proof}}
As detailed in \S\ref{sec-general-properties}, it follows directly from Theorem~\ref{thm-matrix-equation} that  $y$ as defined in (\ref{global-solution}) is the unique standard solution to (\ref{riccati-2}), provided $y_0$ does not belong to the singular curve $\Sigma_M$, where $M=\psi(f,g,h)$.  

Suppose on the other hand that $y_0\in\Sigma_M$, and that $y$ is defined by (\ref{global-solution}).  We claim that $\ster y\in\aclocgen(\real,S^2)$.  Justification of this claim takes some work.  To begin, we make some basic observations concerning stereographic projection, and set up some notation.  Absolute continuity of functions $u:\real\rightarrow S^2$ will be based on the standard euclidean $\real^3$ metric on $S^2\subset\real^3$.  This is equivalent to the metric of geodesic distance on $S^2$ as far as absolute continuity is concerned, but is easier to work with.  Under stereographic projection, the mapping on the Riemann sphere
\begin{equation}\nonumber
\iota:\rsphere\rightarrow\rsphere,\qquad\iota(z)=\frac{1}{z}
\end{equation}
is an isometry of $S^2$. Indeed, for every $(x_1,x_2,x_3)\in S^2$,
\begin{equation}\label{ster-isometry}
\ster\iota\ster^{-1}(x_1,x_2,x_3)=(x_1,-x_2,-x_3). 
\end{equation}
Therefore, given $u:\real\rightarrow\rsphere$, $\ster u\in\aclocgen(\real,S^2)$ if and only if $\ster\frac{1}{u}\in\aclocgen(\real,S^2)$. 
Distance in the complex plane bounds stereographic distance, since
\begin{equation}\label{ster-ac}
\forall z_1,z_2\in\complex,\quad \|\ster(z_1)-\ster(z_2)\|=\frac{2|z_1-z_2|}{\sqrt{(1+|z_1|^2)(1+|z_2|^2)}}\leq 2|z_1-z_2|.
\end{equation}
Therefore, if $u\in\acloc$, then automatically $\ster u\in\aclocgen(\real,S^2)$.  Fix notation for intervals as follows.  For any $x\in\real$ and $r>0$, let $B_r(x)$ denote the open interval $(x-r,x+r)$.  We now proceed to the claim at hand, namely that $\ster y\in\aclocgen(\real,S^2)$. 

Referring to (\ref{psi}), set 
\begin{equation}\label{alpha-beta}
\alpha=e^{Pg}\bigl(y_0C_{F,H}+S_{H,F}\bigr),\quad\beta=y_0S_{F,H}+C_{H,F},\quad\mbox{ so that }\quad y=\frac{\alpha}{\beta}. 
\end{equation}
Note that $\alpha$ and $\beta$ cannot be simultaneously 0, since $\det M\equiv 1$ (where $M=\psi(f,g,h)$, as above).  Denote 
\begin{equation}\label{zeros}
Z_\beta=\beta^{-1}(0)\quad\mbox{ and }\quad Z_\alpha=\alpha^{-1}(0).
\end{equation}
Thus $Z_\alpha,Z_\beta\subset\real$ are disjoint, closed sets.  For each $x\in\real\setminus Z_\beta$, let $r_x>0$ be such that 
\begin{equation}\nonumber
B_{3r_x}(x)\subset\real\setminus Z_\beta, 
\end{equation}
and set $U_x=B_{r_x}(x)$.  For each $x\in Z_\beta$, let $r_x>0$ be such that 
\begin{equation}\nonumber
B_{3r_x}(x)\subset\real\setminus Z_\alpha, 
\end{equation}
and set $V_x=B_{r_x}(x)$.  Thus the family of sets
\begin{equation}\nonumber
\mathscr{F}=\left\{\left.U_x\,\right|\,x\in\real\setminus Z_\beta\right\}\cup\left\{\left.V_x\,\right|\,x\in Z_\beta\right\}
\end{equation}
is an open cover of $\real$.  Fix a nonempty compact interval $K\subset\real$.  The interval $K$ has a finite subcover from $\mathscr{F}$,
\begin{equation}\label{subcover}
K\subset U_{x_1}\cup\ldots\cup U_{x_m}\cup V_{x_{m+1}}\cup\ldots\cup V_{x_{m+n}}.
\end{equation}
Referring to this finite subcover (\ref{subcover}), set 
\begin{equation}\label{script-UV}
\begin{split}
\mathcal{U}_1&=\bigcup_{j=1}^mU_{x_j},\quad\mathcal{U}_2=\bigcup_{j=1}^m\overline{B_{2r_{x_j}}(x_j)},\\
\mathcal{V}_1&=\bigcup_{j=1}^nV_{x_{m+j}},\quad\mathcal{V}_2=\bigcup_{j=1}^n\overline{B_{2r_{x_{m+j}}}(x_{m+j})}.
\end{split}
\end{equation}
Observe that by construction: $\mathcal{U}_1$ and $\mathcal{V}_1$ are open; $K\subset \mathcal{U}_1\cup\mathcal{V}_1$; $\mathcal{U}_2$ is a compact subset of $\real\setminus Z_\beta$ that contains $\mathcal{U}_1$; $\mathcal{V}_2$ is a compact subset of $\real\setminus Z_\alpha$ that contains $\mathcal{V}_1$.  

Next observe that $y=\alpha/\beta$ is absolutely continuous on $\mathcal{U}_2$, since both $\alpha$ and $\beta$ are locally absolutely continuous, and $\beta$ is bounded away from zero on $\mathcal{U}_2$.  Similarly, $1/y$ is absolutely continuous on $\mathcal{V}_2$, since $\alpha$ is bounded away from zero on $\mathcal{V}_2$.  Therefore by (\ref{ster-ac}), $\ster y$ is absolutely continuous on $\mathcal{U}_2$, and $\ster\frac{1}{y}$ is absolutely continuous on $\mathcal{V}_2$.  As a consequence of the isometry (\ref{ster-isometry}), absolute continuity of $\ster\frac{1}{y}$ is equivalent to that of $\ster y$.  Thus $\ster y$ is absolutely continuous on $\mathcal{V}_2$.  

To see that $\ster y$ is absolutely continuous on $K$, let $\varepsilon>0$ be given.  By absolute continuity of $\ster y$ on $\mathcal{U}_2$, there exists $\delta_{\mathcal{U}}>0$ such that for any finite collection of disjoint intervals $(a_j,b_j)$ in $\mathcal{U}_2\cap K$ of total length at most $\delta_{\mathcal{V}}$,
\begin{equation}\label{epsilon-bound-1}
\sum_{j}\|\ster y(a_j)-\ster y(b_j)\|<\varepsilon/2. 
\end{equation}
And, by absolute
continuity of $\ster y$ on $\mathcal{V}_2$, there exists $\delta_{\mathcal{V}}>0$ such that for any finite collection of disjoint intervals $(a_j,b_j)$ in $\mathcal{V}_2\cap K$ of total length at most $\delta_{\mathcal{V}}$, (\ref{epsilon-bound-1}) holds.  Set 
\begin{equation}\label{delta}
\delta=\min\left\{\delta_{\mathcal{U}},\delta_{\mathcal{V}},\min\{r_{x_j}\,|\,1\leq j\leq m+n\}\right\}.
\end{equation}
Let $(a_j,b_j)$ be a finite collection of disjoint intervals in $K$ of total length at most $\delta$.  Suppose a given interval $(a_j,b_j)$ has nonempty intersection with $\mathcal{U}_1$. Then it intersects some $U_{x_k}=B_{r_{x_k}}(x_k)$, where $1\leq k\leq m$.  Since $b_j-a_j<r_{x_k}$, it follows that $(a_j,b_j)$ is contained in $\overline{B_{2r_{x_k}}(x_k)}$.  Thus $(a_j,b_j)$ is contained in $\mathcal{U}_2$.  Similarly, if $(a_j,b_j)$ has nonempty intersection with $\mathcal{V}_1$, it is contained in $\mathcal{V}_2$.  Since $K\subset \mathcal{U}_1\cup\mathcal{V}_1$, every interval $(a_j,b_j)$ in the collection is a subset either of $\mathcal{U}_2$ or $\mathcal{V}_2$.  By choice of $\delta_{\mathcal{U}}$ and $\delta_{\mathcal{V}}$, absolute continuity of $\ster y$ on each of $\mathcal{U}_2$ and $\mathcal{V}_2$ implies 
\begin{equation}\nonumber
\sum_j\|\ster y(a_j)-\ster y(b_j)\|<\frac{\varepsilon}{2}+\frac{\varepsilon}{2}=\varepsilon. 
\end{equation}
This proves that $\ster y$ is absolutely continuous on $K$.  Since $K$ was arbitrary, it follows that $\ster y\in\aclocgen(\real,S^2)$, as claimed. 

Next we claim that if $f^{-1}(0)$ has measure zero, then for every $y_0\in\Sigma_M$, equation (\ref{riccati-2}) holds almost everywhere.  The basic computation underpinning Theorem~\ref{thm-matrix-equation}, i.e., verification that $M=\psi(f,g,h)$ satisfies equation (\ref{matrix-equation-2}), yields that (\ref{riccati-2}) holds almost everywhere on $\real\setminus Z_\beta$. The main point is that, since $y=\alpha/\beta$ is a ratio of locally absolutely continuous functions, $y^\prime(x)$ is well defined at almost every $x$ for which $\beta(x)\neq 0$.  The problem lies with $x\in Z_\beta$, which could in principle have positive measure.  It turns out that if $f^{-1}(0)$ has measure zero, then so does $Z_\beta$, in which case (\ref{riccati-2}) holds almost everywhere. To see this, let $x\in Z_\beta\setminus f^{-1}(0)$.  Then $\alpha(x)\neq 0$, and by Lemma~\ref{lem-differential-formula},
\begin{equation}\nonumber
\beta^\prime(x)=F(x)(y_0C_{F,H}(x)+S_{H,F}(x))=-f(x)\alpha(x)\neq0. 
\end{equation}
Thus for almost every $x\in Z_\beta$, $\beta^\prime(x)$ exists and is nonzero, and so $\beta\neq0$ in a punctured neighbourhood of $x$. In other words, almost every $x\in Z_\beta$ is an isolated point of $Z_\beta$.  It follows that $Z_\beta$ has measure zero, as claimed. 

Lastly, we prove uniqueness of projective solutions.  
Any standard solution $y_1$ to (\ref{riccati-2}) is uniquely determined by its value at any given $x_0\in\real$.  Given a standard solution $y_1$ and a projective solution $y_2$, either they coincide, or they differ at every $x\in\real$.  For, if they coincide at a point, they necessarily coincide in the maximal interval of existence of $y_2$ about that point. If this maximal interval of existence is not $\real$, then $y_2$ tends to $\infty$ at one of the endpoints, which is precluded by coincidence with $y_1$.   

Since for any fixed $(f,g,h)\in\X$ the singular curve $\Sigma_{\psi(f,g,h)}$ has measure zero, there exist three distinct standard solutions $y_1,y_2,y_3$ to (\ref{riccati-2}). 
Suppose $y_4$ is a projective solution to (\ref{riccati-2}).  If for some $j\in\{1,2,3\}$, and some $x_0\in\real$, $y_4(x_0)=y_j(x_0)$ then $y_4=y_j$ is a standard solution.  Suppose therefore that for every $x\in\real$ and every $j\in\{1,2,3\}$,  $y_4(x)\neq y_j(x)$. Set 
\begin{equation}\nonumber
A=\frac{y_4-y_1}{y_4-y_2}\quad\mbox{ and }\quad B=\frac{y_3-y_1}{y_3-y_2}.
\end{equation}
Then for every $x\in\real$, $A(x)$ and $B(x)$ are complex numbers distinct from 0. $B(x)$ is everywhere distinct from 1, and $A(x)=1$ only where $y_4(x)=\infty$. Since each $y_j$ satisfies (\ref{riccati-2}) almost everywhere, for the same $f,g,h$ but different initial $y_0$, it is straightforward to verify that 
\begin{equation}\nonumber
\left(\frac{A}{B}\right)^\prime=\frac{BA^\prime-AB^\prime}{B^2}=0\quad\mbox{ almost everywhere, }
\end{equation}
and hence that the cross ratio $A/B\equiv \gamma $ for some constant $\gamma \neq 0$.  The function $y_4$ therefore satisfies the formula
\begin{equation}\label{y4-formula}
y_4=\frac{y_1-\gamma y_2B}{1-\gamma B}.
\end{equation}
The right-hand side of (\ref{y4-formula}) is finite at points $x\in\real$ at which $\gamma B(x)\neq 1$, or equivalently, at which $A(x)\neq 1$.  But
$A(x)=1$ if and only if $y_4(x)=\infty$.  Since $y_4$ is assumed to be a projective solution, $y_4(x)=\infty$ on a set of measure zero.  Thus the formula (\ref{y4-formula}) determines $y_4$ almost everywhere, and hence everywhere, by continuity of $\ster y_4$.  Moreover, if $y_4(0)=y_0\in\complex$, this value determines $\gamma$ and hence the value of $y_4$ everywhere.  Thus every projective solution to (\ref{riccati-2}) is uniquely determined by its initial value. 

This completes the proof of Theorem~\ref{thm-riccati-solution}.

\section{Examples and applications\label{sec-examples}}

\subsection{A special case considered by Abel\label{sec-abel}}
Theorem~\ref{thm-riccati-solution} provides a systematic way to recover formulas previously derived for special cases by ad hoc methods.  As an illustration, consider the case of (\ref{riccati}) considered by Abel \cite[\S{IV}]{Ab:1881}, in which 
\begin{equation}\label{abel-case}
h=\gamma^2fe^{2Pg}\mbox{ for a fixed nonzero constant }\gamma\in\complex.
\end{equation}
According to (\ref{psi}), in this case $F=-e^{Pg}f$ and $H=e^{-Pg}h=-\gamma^2F$.  It follows from the defining formulas (\ref{Aj}) and (\ref{cs}), and equation (\ref{diagonal-values}), that 
\begin{equation}\nonumber
C_{F,H}=C_{F,-\gamma^2F}=C_{i\gamma F,i\gamma F}=\cosh(i\gamma PF)=\cos(\gamma PF)\quad\mbox{ and }\quad C_{H,F}=C_{F,H}.
\end{equation}
Also,
\begin{equation}\nonumber
S_{F,H}=S_{F,-\gamma^2F}=\frac{1}{i\gamma}S_{i\gamma F,i\gamma F}=\frac{1}{i\gamma}\sinh(i\gamma PF)=\frac{1}{\gamma}\sin(\gamma PF)
\end{equation}
and
\begin{equation}\nonumber
S_{H,F}=S_{-\gamma^2F,F}={i\gamma}S_{i\gamma F,i\gamma F}={i\gamma}\sinh(i\gamma PF)=-{\gamma}\sin(\gamma PF).
\end{equation}
The formula (\ref{y-full}) therefore yields
\begin{equation}\label{abel-solution}
\begin{split}
y&=e^{Pg}\frac{y_0\cos(\gamma PF)-\gamma\sin(\gamma PF)}{\frac{y_0}{\gamma}\sin(\gamma PF)+\cos(\gamma PF)}\\
&=\gamma e^{Pg}\frac{\frac{y_0}{\gamma}+\tan\bigl(\gamma P(e^{Pg}f)\bigr)}{1-\frac{y_0}{\gamma}\tan\bigl(\gamma P(e^{Pg}f)\bigr)}\\
&=\gamma e^{Pg}\tan\bigl(\gamma P(e^{Pg}f)+\tan^{-1}(y_0/\gamma)\bigr). 
\end{split}
\end{equation}
This generalizes to arbitrary initial values the solution given by Abel.

\subsection{The Airy equation\label{sec-Airy}}
Applied to the Airy equation
\begin{equation}\label{airy-equation}
u^{\prime\prime}-xu=0,
\end{equation}
Theorem~\ref{thm-linear} yields
\begin{equation}\label{airy-solution}
u=u(0)C_{x,1}+u^{\prime}(0)S_{1,x}.
\end{equation}
That each of the functions
\begin{equation}\label{airy-functions}
u_1=C_{x,1}\quad\mbox{ and }\quad u_2=S_{1,x}
\end{equation}
satisfies the Airy equation (\ref{airy-equation}) is easily confirmed directly using Lemma~\ref{lem-differential-formula}. 
Definition~\ref{defn-bivariate-exponential} provides explicit formulas
\begin{equation}\label{airy-formula-1}
u_1(x)=\left\{
\begin{array}{cc}
1+\sum_{j=1}^\infty\thickspace\displaystyle\int\limits_{0<s_1<\cdots<s_{2j}<x}s_1s_3\cdots s_{2j-1}\,ds_1ds_2\cdots ds_{2j}&\mbox{ if }x\geq0\\
1+\sum_{j=1}^\infty\thickspace\displaystyle\int\limits_{x<s_{2j}<\cdots<s_1<0}s_1s_3\cdots s_{2j-1}\,ds_1ds_2\cdots ds_{2j}&\mbox{ if }x\leq0
\end{array}
\right.
\end{equation}
\begin{equation}\label{airy-formula-2}
u_2(x)=
\left\{\begin{array}{cc}
\sum_{j=1}^\infty\thickspace\displaystyle\int\limits_{0<s_1<\cdots<s_{2j-1}<x}s_2s_4\cdots s_{2j-2}\,ds_1ds_2\cdots ds_{2j-1}&\mbox{ if }x\geq0\\
-\sum_{j=1}^\infty\thickspace\displaystyle\int\limits_{x<s_{2j-1}<\cdots<s_1<0}s_2s_4\cdots s_{2j-2}\,ds_1ds_2\cdots ds_{2j-1}&\mbox{ if }x\leq0
\end{array}
\right.
\end{equation}
The formulas (\ref{airy-formula-1}) and (\ref{airy-formula-2}) for Airy functions appear to be new. 

The values of the standard Airy functions and their derivatives at the origin are
\begin{equation}\label{standard-airy-values}
{\rm Ai}(0)=\frac{1}{3^{\frac{2}{3}}\Gamma(\frac{2}{3})},\quad  {\rm Ai}^\prime(0)=-\frac{1}{3^{\frac{1}{3}}\Gamma(\frac{1}{3})},\quad  
{\rm Bi}(0)=\frac{1}{3^{\frac{1}{6}}\Gamma(\frac{2}{3})},\quad  {\rm Bi}^\prime(0)=\frac{3^{\frac{1}{6}}}{\Gamma(\frac{1}{3})}.
\end{equation}
See \cite[p.~91]{Wo:2001}. It follows that Ai and Bi are given 
in terms of $u_1$ and $u_2$ 
by 
\begin{equation}\label{standard-airy}
{\rm Ai}=\frac{u_1}{3^{\frac{2}{3}}\Gamma(\frac{2}{3})}-\frac{u_2}{3^{\frac{1}{3}}\Gamma(\frac{1}{3})}\quad\mbox{ and }\quad
{\rm Bi}=\frac{u_1}{3^{\frac{1}{6}}\Gamma(\frac{2}{3})}+\frac{3^{\frac{1}{6}}u_2}{\Gamma(\frac{1}{3})}.
\end{equation}

\subsection{The Schr\"odinger equation\label{sec-schroedinger}}
Given $q\in\lloc$ and $\lambda\in\complex$, application of Theorem~\ref{thm-linear} to the one-dimensional Schr\"odinger equation
\begin{equation}\label{schroedinger}
-y^{\prime\prime}+qy=\lambda y
\end{equation}
yields 
\begin{equation}\label{schoedinger-solution}
y=y(0)C_{q-\lambda,1}+y^\prime(0)S_{1,q-\lambda}. 
\end{equation}
Definition~\ref{defn-bivariate-exponential} renders this solution explicit.  
The first explicit solution was given only recently, by a much more complicated derivation \cite{Gi:SIMA2024}. 

\subsection{The Miura transform\label{sec-miura}}

For real valued $q\in\lloc$ and $\alpha\in\acloc$, inversion of the Miura transform
\begin{equation}\label{miura-transform}
\alpha\mapsto q=\alpha^2+\alpha^\prime, 
\end{equation}
is equivalent to determination of real valued solutions to the Riccati equation 
\begin{equation}\label{miura-riccati}
y^\prime=q-y^2.
\end{equation}
Theorem~\ref{thm-riccati-solution} tells us how to do this. Here $(f,g,h)=(-1,0,q)$, and
\begin{equation}\label{miura-case-psi}
M=\psi(-1,0,q)=\begin{pmatrix}C_{1,q}&S_{q,1}\\ S_{1,q}&C_{q,1}\end{pmatrix}.
\end{equation}
So, the real valued solutions to (\ref{miura-riccati}) are precisely the functions
\begin{equation}\label{miura-solution}
y=\frac{y_0C_{1,q}+S_{q,1}}{y_0S_{1,q}+C_{q,1}}\quad\mbox{ where }\quad y_0\in\real,
\end{equation}
each of which is guaranteed to be at least a projective solution, 
since $f=-1$ has no zeros.  Thus, given a real valued function $q\in\lloc$, the number $y_0=\alpha(0)\in\real$ determines by (\ref{miura-solution}) a unique function $\alpha=y$ such that $q=\alpha^2+\alpha^\prime$ almost everywhere and $\ster\alpha\in\aclocgen(\real,S^2)$.  
It turns out that the set of initial values $y_0\in\real$ giving rise to a locally absolutely continuous solution to (\ref{miura-riccati}) is connected, as follows. 
\begin{prop}\label{prop-interval}
Let $q\in\lloc$ be real valued, and set $M=\psi(-1,0,q)$ as per (\ref{psi}).  Then $\Sigma_M\subset\real$.  If $\Sigma_M\neq\real$, then $\real\setminus\Sigma_M$ is a connected interval. 
\end{prop}
\begin{pf}
The singular curve associated to $M=\psi(-1,0,q)$ is
\begin{equation}\label{sigma-form}
\Sigma_M=\left\{ \sigma(x)\,\left|\,x\in\real\,\&\,S_{1,q}(x)\neq0\right.\right\}\quad\mbox{ where }\quad\sigma=-C_{q,1}/S_{1,q}.
\end{equation}
Since $q$ is real valued, it follows immediately that $\Sigma_M\subset\real$. 

The set $Z_q=S_{1,q}^{-1}(0)$ is a closed subset of the real line containing 0.  Thus every $x_0\not\in Z_q$ belongs to a maximal interval $(\alpha,\beta)$ disjoint from $Z_q$, at least one of whose endpoints is finite; i.e., $\alpha\neq-\infty$ or $\beta\neq\infty$.  If $\alpha\neq-\infty$, then $S_{1,q}(\alpha)=0$.  Since $\det M\equiv 1$, it follows that $C_{1,q}(\alpha)\neq 0$, and so 
\begin{equation}\nonumber
\lim_{x\rightarrow\alpha+}\sigma(x)=\pm\infty. 
\end{equation}
Thus the interval $\sigma\bigl((\alpha,\beta)\bigr)$ has $\pm\infty$ as an endpoint, and therefore its complement 
\begin{equation}\nonumber
S_{x_0}=\real\setminus\sigma\bigl((\alpha,\beta)\bigr)
\end{equation}
is connected, or empty.  By the same argument, the set $S_{x_0}$ is connected if $\beta\neq\infty$.  The complement $\real\setminus\Sigma_M$ of the singular curve is the intersection over $x_0\in\real\setminus Z_q$ of connected sets $S_{x_0}$.  Therefore $\real\setminus\Sigma_M$ is also connected, provided it is nonempty. 
\end{pf}

Initial values $y_0\in\real\setminus\Sigma_M$ correspond via (\ref{miura-riccati}) to continuous pre-images of $q$ by the Miura transform. All other initial values generate discontinuous pre-images.  
Kappeler et al.~obtain a result akin to Proposition~\ref{prop-interval} by completely different arguments in \cite[Thm.~1.3]{KaPeShTo:2005}.  The explicit characterization of $\Sigma_M$ given here provides an alternate, more concrete perspective.

\section{Remarks}\label{sec-remarks}

\subsection{On the bivariate exponential operator\label{sec-bivariate-exponential-2}}

The close relation of the bivariate exponential to the familiar exponential primitive operator, the simplicity of its formulation, and the evident utility its differential formulas as expressed in Lemma~\ref{lem-differential-formula}, all suggest that the bivariate exponential operator is a fundamental object that ought to have been known long ago.  Yet it appears to be new.  We were led to it by a rather circuitous route, via the impedance form of the Schr\"odinger equation.  As described in \cite{Gi:SIMA2024}, there is a much more restricted precursor to the bivariate exponential, called the harmonic exponential operator, that satisfies a differential equation obtained as a renormalized continuum limit of the three-term recurrence formula for orthogonal polynomials on the circle.  
However, this somewhat arcane provenance is a historical accident rather than a reflection of mathematical necessity; the bivariate exponential of Definition~\ref{defn-bivariate-exponential} could just as easily have been written down from scratch.

\subsection{Conjugacy classes of Riccati equations\label{sec-conjugacy-classes}}
Established theory for the Riccati equation includes a miscellany of
ways to determine a new solution, given a particular solution.  The formulation of Theorem~\ref{thm-matrix-equation} suggests a natural organizing principle for interrelated solutions to (\ref{riccati}) by conjugacy class, as follows.

Note that conjugation by an arbitrary invertible matrix $J\in\GL(2,\complex)$ leaves both $\mathfrak{sl}(2,\complex)$ and $\SL(2,\complex)$ invariant. Suppose $J\in\GL(2,\complex)$, $\mathfrak{m}_1\in\C$ and $\mathfrak{m}_2=J\mathfrak{m}_1J^{-1}$.  Denote 
\begin{equation}\nonumber
\mathfrak{m}_j=\begin{pmatrix}\frac{1}{2}g_j&h_j\\ -f_j&-\frac{1}{2}g_j\end{pmatrix}\qquad (j=1,2),
\end{equation}
so that $(f_1,g_1,h_1)$ and $(f_2,g_2,h_2)$ are related according to the conjugation relation of $\mathfrak{m}_1$ and $\mathfrak{m}_2$. 
Suppose $y_1$ satisfies the Riccati equation 
\begin{equation}\nonumber
y^\prime=f_1y^2+g_1y+h_1,\quad y(0)=y_0.
\end{equation}
Then it follows from Theorems~\ref{thm-matrix-equation} and \ref{thm-riccati-solution} that the solution $y_2$ to the conjugate equation
\begin{equation}\nonumber
y^\prime=f_2y^2+g_2y+h_2,\quad y(0)=\varphi_{J}(y_0)
\end{equation}
is given by the formula 
\begin{equation}\label{J-formula}
y_2(x)=\varphi_J(y_1(x))\qquad (x\in\real). 
\end{equation}
Thus, referring to (\ref{spaces}), pairs of elements of $\X$ whose corresponding matrices in $\C$ are conjugate are coefficients of Riccati equations whose solutions are related by a constant linear fractional transformation.  It follows that if a given Riccati equation is soluble by an elementary combination of quadratures, then the same is true of its entire conjugacy class.  Interestingly, no two distinct equations of the form (\ref{miura-riccati}) are conjugate to one another. 

In the particular case $J=\begin{pmatrix}0&1\\ 1&0\end{pmatrix}$, we have 
\begin{equation}\nonumber
(f_2,g_2,h_2)=(-h_1,-g_1,-f_1)\quad\mbox{ and }\quad y_2=\varphi_J(y_1)=\frac{1}{y_1}.
\end{equation}
Thus application of the involution 
\begin{equation}\nonumber
(f,g,h)\mapsto(-h,-g,-f)
\end{equation}
to coefficients of the Riccati equation (\ref{riccati}) corresponds to transformation of the solution by the Riemann sphere isometry $y\mapsto 1/y$.  

\subsection{Weak solutions and lower regularity of the coefficients\label{sec-weak-solutions}}
Bijectivity of the Riccati correspondence affirms $\X$ as the natural domain for coefficients $(f,g,h)$ occurring in (\ref{riccati}).  Historically, coefficients are typically assumed to be continuous, and in many cases real valued.  The latter more restrictive assumptions are in part a reflection of the Riccati equation's antiquity, pre-dating complex analysis and the Lebesgue theory of measure and integration;
the assumption that coefficients be real is also related to physical applications.  In addition, it is commonly assumed that the coefficient $f$ of $y^2$ is nonzero, since (\ref{riccati}) reduces to a linear equation on intervals where $f(x)=0$.  From this point of view, the hypothesis in Theorem~\ref{thm-riccati-solution} that $f^{-1}(0)$ have measure zero can be seen as a condition of nonlinearity of (\ref{riccati}). 

It is possible to study the Riccati equation with coefficients of lower regularity than $\lloc$, in which case one has to consider some form of weak solution.  This is relevant both in the context of the Schr\"odinger equation and the Miura map, where it is of interest to allow $q\in H^{-1}_{\mbox{\tiny loc}}(\real)$ \cite{KaPeShTo:2005} or even $q\in H^{-2}_{\mbox{\tiny loc}}$ \cite{Gi:SIMA2024}.  
For the Schr\"odinger equation in impedance form, approximation by piecewise constant impedance corresponds to the pre-image of $q$ by the Miura map having the form of a sum of Dirac delta functions. This in turn has interesting connections to the geometry of $\aut\rsphere$ \cite[\S2\&\S6]{Gi:JAT2023}. 
These considerations call for further investigation, beyond the scope of the present paper. 

\section*{Acknowledgements}
Research was partially supported by NSERC grant DG RGPIN-2022-04547.

%


\end{document}